\documentclass[12pt]{amsart}
\usepackage{amssymb,latexsym}
\usepackage{pdfsync}

\newdimen\AAdi%
\newbox\AAbo%
\def\AAk#1#2{\s_etbox\AAbo=\hbox{#2}\AAdi=\wd\AAbo\kern#1\AAdi{}}%
\def\AAr#1#2#3{\s_etbox\AAbo=\hbox{#2}\AAdi=\ht\AAbo\raise#1\AAdi\hbox{#3}}%
\font\tenmsb=msbm10 at 12pt
\font\sevenmsb=msbm7 at 8pt
\font\fivemsb=msbm5 at 6pt
\newfam\msbfam
\textfont\msbfam=\tenmsb
\scriptfont\msbfam=\sevenmsb
\scriptscriptfont\msbfam=\fivemsb
\def\Bbb#1{{\tenmsb\fam\msbfam#1}}

\newcommand{\beq}{\begin{equation}}
\newcommand{\eeq}{\end{equation}}
\newcommand{\beqr}{\begin{eqnarray}}
\newcommand{\eeqr}{\end{eqnarray}}
\newcommand{\ba}{\begin{array}}
\newcommand{\ea}{\end{array}}

\begin{document}

\newtheorem{thm}{Theorem}
\newtheorem{lem}{Lemma}
\newtheorem{cor}{Corollary}
\newtheorem{rem}{Remark}
\newtheorem{pro}{Proposition}
\newtheorem{defi}{Definition}
\newtheorem{eg}{Example}
\newtheorem*{claim}{Claim}
\newtheorem{conj}[thm]{Conjecture}
\newcommand{\noi}{\noindent}
\newcommand{\dis}{\displaystyle}
\newcommand{\mint}{-\!\!\!\!\!\!\int}
\numberwithin{equation}{section}

\def \bx{\hspace{2.5mm}\rule{2.5mm}{2.5mm}}
\def \vs{\vspace*{0.2cm}}
\def\hs{\hspace*{0.6cm}}
\def \ds{\displaystyle}
\def \p{\partial}
\def \O{\Omega}
\def \o{\omega}
\def \b{\beta}
\def \m{\mu}
\def \l{\lambda}
\def\L{\Lambda}
\def \ul{u_\lambda}
\def \D{\Delta}
\def \d{\delta}
\def \k{\kappa}
\def \s{\sigma}
\def \e{\varepsilon}
\def \a{\alpha}
\def \tf{\tilde{f}}
\def\cqfd{%
\mbox{ }%
\nolinebreak%
\hfill%
\rule{2mm} {2mm}%
\medbreak%
\par%
}
\def \pr {\noindent {\it Proof.} }
\def \rmk {\noindent {\it Remark} }
\def \esp {\hspace{4mm}}
\def \dsp {\hspace{2mm}}
\def \ssp {\hspace{1mm}}

\def\la{\langle}\def\ra{\rangle}

\def \u{u_+^{p^*}}
\def \ui{(u_+)^{p^*+1}}
\def \ul{(u^k)_+^{p^*}}
\def \energy{\int_{\R^n}\u }
\def \sk{\s_k}
\def \mo{\mu_k}
\def\cal{\mathcal}
\def \I{{\cal I}}
\def \J{{\cal J}}
\def \K{{\cal K}}
\def \OM{\overline{M}}

\def\n{\nabla}

\def\fk{{{\cal F}}_k}
\def\M1{{{\cal M}}_1}
\def\Fk{{\cal F}_k}
\def\Fl{{\cal F}_l}
\def\FF{\cal F}
\def\Gk{{\Gamma_k^+}}
\def\n{\nabla}
\def\uuu{{\n ^2 u+du\otimes du-\frac {|\n u|^2} 2 g_0+S_{g_0}}}
\def\uuug{{\n ^2 u+du\otimes du-\frac {|\n u|^2} 2 g+S_{g}}}
\def\sku{\sk\left(\uuu\right)}
\def\qed{\cqfd}
\def\vvv{{\frac{\n ^2 v} v -\frac {|\n v|^2} {2v^2} g_0+S_{g_0}}}
\def\vvs{{\frac{\n ^2 \tilde v} {\tilde v}
 -\frac {|\n \tilde v|^2} {2\tilde v^2} g_{S^n}+S_{g_{S^n}}}}
\def\skv{\sk\left(\vvv\right)}
\def\tr{\hbox{tr}}
\def\pO{\partial \Omega}
\def\dist{\hbox{dist}}
\def\RR{\Bbb R}\def\R{\Bbb R}
\def\C{\Bbb C}
\def\B{\Bbb B}
\def\N{\Bbb N}
\def\Q{\Bbb Q}
\def\Z{\Bbb Z}
\def\PP{\Bbb P}
\def\EE{\Bbb E}
\def\F{\Bbb F}
\def\G{\Bbb G}
\def\H{\Bbb H}
\def\SS{\Bbb S}\def\S{\Bbb S}

\def\div{\hbox{div}\,}

\def\lcf{{locally conformally flat} }

\def\circledwedge{\setbox0=\hbox{$\bigcirc$}\relax \mathbin {\hbox
to0pt{\raise.5pt\hbox to\wd0{\hfil $\wedge$\hfil}\hss}\box0 }}

\def\sss{\frac{\s_2}{\s_1}}

\date{\today}
\title[ A rigidity theorem of ancient solutions  ]{ A rigidity theorem of ancient solutions to the mean curvature flow in codimension one}

\author{}

\author[Chen]{Qun Chen}
\address{School of Mathematics and Statistics\\ Wuhan University\\Wuhan 430072,
China
 }
 \email{qunchen@whu.edu.cn}
 
\author[Qiu]{Hongbing Qiu$^*$}
\address{School of Mathematics and Statistics\\ Wuhan University\\Wuhan 430072,
China
 }
 \email{hbqiu@whu.edu.cn}

 \renewcommand{\thefootnote}{\fnsymbol{footnote}}

\footnotetext{\hspace*{-5mm} \begin{tabular}{@{}r@{}p{13.4cm}@{}}

 $^*$ & Corresponding author \\

\end{tabular}}

\renewcommand{\thefootnote}{\arabic{footnote}}
 

\begin{abstract}

By carrying out a  point-wise estimate for the second fundamental form, we prove a rigidity theorem of complete noncompact ancient solutions to the mean curvature flow in codimension one. Moreover, we derive an optimal growth condition.

\vskip12pt



\end{abstract}
\maketitle
\section{Introduction}

Let $F: M^{n} \rightarrow \mathbb{R}^{m+n}$ be an isometric immersion from an $n$-dimensional oriented Riemannian manifold $M$ to the Euclidean space $\mathbb{R}^{m+n}$.  
The mean curvature flow (MCF) in Euclidean space is a one-parameter family of immersions $F_t= F(\cdot, t): M^n \rightarrow \mathbb{R}^{m+n}$ with the corresponding images $M_t=F_t(M)$ such that
\begin{equation}\label{eqn-MCF}
\begin{cases}\aligned
\frac{\partial}{\partial t}F(x,t)=& H(x,t), x\in M,\\
F(x,0)=&F(x),
\endaligned
\end{cases}
\end{equation}
is satisfied, where $H(x, t)$ is the mean curvature vector of $M_t$ at $F(x, t)$ in $\mathbb{R}^{m+n}$.

A solution of the mean curvature flow, is called ancient if it exists on a time interval of the form $(-\infty, T_0)$ with $T_0 < \infty$. It 
models the asymptotic profile of the MCF near a singularity (c.f. \cite{HuiSin199, HuiSin299}). Therefore, determining the possible shapes of ancient solutions is a central topic in the singularity analysis. In the past decade, there has been important progress on the classification of compact ancient solutions to the MCF, for instance, see \cite{ADS19, ADS20, BIS17, BL16, BLT21, CM19, CM22, DHS10, HH16, HS15, Lan17, Lan17a, LL20, LN21, LXZ21, SS19, Wan11a, Whi03}.  

For the noncompact ancient solutions, Brendle-Choi  \cite{BC19} classified the noncompact ancient non-collapsed solution to the MCF in $\mathbb R^3$, and more generally in high dimensional Euclidean space $\mathbb R^{n+1}$ under an additional assumption of uniform two-convexity  \cite{BC21}. Recently, Haslhofer and coauthors \cite{CH21, DH21, DH23, DH24} gave the classification of the ancient noncollapsed flows in $\mathbb R^4$.  In \cite{Kun16}, Kunikawa showed that there are no nontrivial ancient solutions with bounded slope and bounded mean curvature in codimension one.  
 On the one hand, this result was generalized  to higher codimension by the  author and Xin \cite{QX23}, that is, any ancient solution in higher codimension has to be an affine subspace provided the slope is smaller than 3. 
See also the related work \cite{GXZ23}. On the other hand, the condition on the mean curvature mentioned above in the case of codimension one could be removed (c.f. \cite{Qiu22}). 

In this note, we continue to study complete noncompact ancient solution to the MCF in codimension one.  With the aid of an appropriate auxiliary function, we derive the following local estimate for the second fundamental form by using the technique of \cite{SZ} and \cite{KS21}. 

\begin{thm} \label{thm-main1}

Let $F: M^n \times [-T, 0] \to \mathbb R^{n+1} $be a solution to the mean curvature flow.  Let $\rho$ be the distance function on $\mathbb S^n$ from some point $p_0$. Assume that the image of the initial submanifold under the Gauss map $\gamma$ is contained in an open hemisphere of $\mathbb S^n$.  
Then we have
\begin{equation}\label{eqn-Pa01}
\aligned
 \sup_{D_{\frac{R}{2}, \frac{T}{2}}(o)} \frac{|B|}{b - \varphi\circ \gamma} 
\leq & C \left( \frac{1}{R}\cdot \left( \sup_{D_{R, T}(o)} \left( \frac{\pi}{2} - \rho\circ \gamma \right)^{-1} \right.\right. \\
& \left.\left.+  \sup_{D_{R, T}(o)} \left( \frac{\pi}{2} - \rho\circ \gamma \right)^{-2} \right) 
+\frac{1}{\sqrt T}\cdot  \sup_{D_{R, T}(o)} \left( \frac{\pi}{2} - \rho\circ \gamma \right)^{-1} \right),
\endaligned
\end{equation} 
 where $b$ and $\varphi$ are defined by (\ref{eqn-b}), and $C>0$ is a constant independent of $R$ and $T$.

\end{thm}

The above estimate (\ref{eqn-Pa01}) leads to a rigidity theorem for ancient solutions to the MCF.

\begin{thm}\label{thm-main2}

Let $F: M^n \times (-\infty, 0] \to \mathbb R^{n+1} $be a complete ancient solution to the mean curvature flow. 
Let $\rho$ be the distance function on $\mathbb S^n$ from some point $p_0$. If, as $t \to -\infty$, the image of the Gauss map $\gamma: M^n \times (-\infty, 0] \to \mathbb S^n$ is contained in an open hemisphere of $\mathbb S^n$ and 
\[
\left( \frac{\pi}{2} - \rho\circ \gamma(x, t) \right)^{-1} = o(\sqrt{|F| + \sqrt{|t|}})
\]
near infinity.  Then $M_t$ has to be affine linear for any $t\in (-\infty, 0]$. 

\end{thm}

\begin{rem}

The author \cite{Qiu22} showed that the same conclusion holds under the condition that the image of the initial submanifold under the Gauss map is contained in a regular ball of $\mathbb S^n$. That is, $\rho\circ \gamma$ is bounded. Hence the condition in the above Theorem \ref{thm-main2} is weaker than the one in \cite{Qiu22}. See also the work by Kunikawa \cite{Kun16}.

\end{rem}

\begin{rem}

Consider the grim reaper: $X=(-\ln \cos x, x, y), |x|<\frac{\pi}{2}, y\in \mathbb R$, by direct computation, we know that its unit normal is $\nu=(\cos x, \sin x, 0)$, namely, the image of  the grim reaper under the Gauss map is a great circle of $\mathbb S^2$. Hence the condition that the image of the Gauss map is contained in an open hemisphere in Theorem \ref{thm-main2} is necessary. 

\end{rem}

\begin{rem}

There exists a well-known entire rotationally symmetric graphic translating soliton $M^2:= \{ F= (x, u(x)) | x\in \mathbb R^2 \} \subset \mathbb R^3$ which grows quadratically at infinity, namely, $|u(x)| \sim C|x|^2$ for some constant $C>0$ ( see for example Lemma 2.2 in \cite{CSS07}), and $\left( \cos(\rho\circ \gamma)\right)^{-1} = \la \nu, (0, 0, 1) \ra^{-1} = O(|F|^{\frac{1}{2}})$, where $|F|=\sqrt{|x|^2 + |u(x)|^2}$ and $\nu=\frac{(-Du, 1)}{\sqrt{1+|Du|^2}}$ is the unit normal of $M^2$ in $\mathbb R^3$. That is, $\left( \frac{\pi}{2} - \rho\circ \gamma \right)^{-1} = O(|F|^{\frac{1}{2}})$. Notice that a translating soliton is a special class of ancient solutions,  which can be viewed as a hypersurface at some time satisfying (\ref{eqn-T111}), thus the growth condition in Theorem \ref{thm-main2} is optimal.

\end{rem}

\bigskip

Recall that $M^n$ is said to be a translating soliton in
$\mathbb{R}^{n+1}$ if it
 satisfies
\begin{equation}\label{eqn-T111}
H= V_{0}^N,
\end{equation}
where  $V_0$ is a fixed vector in $\mathbb{R}^{n+1}$ with unit
length and $V_{0}^N$ denotes the orthogonal projection of $V_0$ onto
the normal bundle of $M^n$.

It is well known that the translating soliton is an important example of the ancient solution to the mean curvature flow.  
In fact, translating solitons do not change their shape under the mean curvature flow, and hence we only have to consider a time slice (a hypersurface at some time) which satisfies (\ref{eqn-T111}). 
Therefore Theorem \ref{thm-main2}  implies the following result of translating solitons, which is  a time independent version of the above Theorem \ref{thm-main2}.

\begin{cor}\label{cor1}

Let $F: M^n \to \mathbb R^{n+1}$ be a complete $n$-dimensional translating soliton  in $\mathbb R^{n+1}$. If  the image of $M^n$ under the Gauss map $\gamma$ is contained in an open hemisphere of $\mathbb S^n$ and 
\[
\left( \frac{\pi}{2} - \rho\circ \gamma \right)^{-1} =o(|F|^{\frac{1}{2}}),
\]
where  $\rho$ is the distance function on $\mathbb S^n$. Then $M^n$ has to be an affine subspace.

\end{cor}

\begin{rem}

The above Corollary \ref{cor1}  is essentially due to the work of Kunikawa \cite{Kun15}.

\end{rem}

\vskip24pt

\section{Preliminaries}

The second fundamental form $B$ of $M^{n}$ in $\mathbb{R}^{n+1}$ is defined by
\[
B_{UW}:= (\overline{\n}_U W)^N
\]
for $U, W \in \Gamma(TM^n)$. We use the notation $( \cdot )^T$ and $(
\cdot )^N$ for the orthogonal  projections into the tangent bundle
$TM^n$ and the normal bundle $NM^n$, respectively. For $\nu \in
\Gamma(NM^n)$ we define the shape operator $A^\nu: TM^n \rightarrow TM^n$
by
\[
A^\nu (U):= - (\overline{\n}_U \nu)^T
\]
Taking the trace of $B$ gives the mean curvature vector $H$ of $M^n$
in $\mathbb{R}^{n+1}$ and
\[
H:= \hbox{trace} (B) = \sum_{i=1}^{n} B_{e_ie_i},
\]
where $\{ e_i \}$ is a local orthonormal frame field of $M^n$.

\bigskip

Let $B_R(o)$ be a Euclidean closed ball of radius $R$ centered at the origin $o \in \mathbb{R}^{n+1}$ and $B_{R, T}(o):= B_R(o) \times [-T, 0]  \subset \mathbb{R}^{n+1}\times (-\infty, +\infty)$. We may consider $\mathcal{M}_T$ as the space-time domain 
\[
\{ (F(p, t), t): p\in M^n, t\in [-T, 0] \} \subset \mathbb{R}^{n+1} \times (-\infty, \infty).
\]
Let $D_{R, T}(o) = \{ (x, t) \in M^n \times [-T, 0]: F(x, t) \in B_R(o) \}.$

\bigskip

Let $\rho$ be the distance function on $\mathbb S^n$ from some point $p_0$. Denote
 \begin{equation}\label{eqn-b}
 \varphi :=1-\cos \rho,  \quad \quad b:=\frac{1}{2}\left( 1+\sup_{D_{R, T}(o)}\varphi\circ \gamma \right).
\end{equation}

\vskip24pt

\section{Proofs of main Theorems }

\vskip12pt

\noindent{\bf Proof of Theorem \ref{thm-main1}} \quad
Direct computation gives us
\begin{equation}\label{eqn-Para00}
{\rm Hess}(\rho) =  \cot\rho\cdot (g - d\rho \otimes d\rho),
\end{equation}
where $g$ is the metric tensor on $\mathbb S^n$. 

Let $\Phi:= 1+\e_0 -\cos \rho$, where $\e_0$ is a fixed positive constant. Then we get
\[
{\rm Hess}(\Phi) = (\cos \rho) g.
\]
By Theorem A in \cite{Wan03}, 
\[
\tau(\gamma) - \p_t \gamma = 0.
\]
It follows that
\begin{equation}\label{eqn-VLap-h101}\aligned
(\p_t-\D) (\Phi\circ \gamma) =& -\sum_{i=1}^m {\rm Hess} (\Phi) (d\gamma(e_i), d\gamma(e_i)) - d\Phi(\tau(\gamma)-\p_t \gamma) \\
=&  -\sum_{i=1}^m {\rm Hess} (\Phi) (d\gamma(e_i), d\gamma(e_i)) = -\cos\rho |d\gamma|^2.
\endaligned
\end{equation}
Hence by the maximum principle we can conclude that the condition $\gamma(M_t) \subset $ an open hemisphere of $\mathbb S^n$ is remained true under the mean curvature flow provided it is valid in the initial submanifold.  

By (\ref{eqn-Para00}),  we obtain
\begin{equation}\label{eqn-Para1phi}\aligned
{\rm Hess}(\varphi) = & \varphi' {\rm Hess} (\rho) + \varphi'' d\rho\otimes d\rho \\
=& \sin\rho{\rm Hess} (\rho) + \cos\rho d\rho\otimes d\rho \\
= & \sin\rho \cot\rho\cdot(g - d\rho \otimes d\rho) + \cos\rho d\rho\otimes d\rho \\
= &( \cos\rho) g.
\endaligned
\end{equation}

Let
\[
f:= \frac{|B|^2}{(b-\varphi\circ \gamma)^2}.
\]
Direct computation gives us
\begin{equation}\label{eqn-Pa1}\aligned
\left( \D -\p_t \right) f = & \frac{\left( \D -\p_t \right) |B|^2}{(b-\varphi\circ \gamma)^2} + \frac{2|B|^2\left(\D -\p_t\right)(\varphi\circ \gamma)}{(b-\varphi\circ \gamma)^3} \\
& + \frac{4\la \n |B|^2, \n(\varphi\circ \gamma) \ra}{(b-\varphi\circ \gamma)^3} + \frac{6|B|^2|\n(\varphi\circ \gamma)|^2}{(b-\varphi\circ \gamma)^4}.
\endaligned
\end{equation}
Since
\begin{equation}\label{eqn-Pa2}\aligned
\frac{\la \n f, \n(\varphi\circ \gamma) \ra}{b-\varphi\circ \gamma} =  \frac{\la \n |B|^2, \n(\varphi\circ \gamma) \ra}{(b-\varphi\circ \gamma)^3} + \frac{2|B|^2|\n(\varphi\circ \gamma)|^2}{(b-\varphi\circ \gamma)^4}.
\endaligned
\end{equation}
From (\ref{eqn-Pa1}) and (\ref{eqn-Pa2}), we derive
\begin{equation}\label{eqn-Pa3}\aligned
\left( \D -\p_t \right) f = & \frac{\left( \D -\p_t \right) |B|^2}{(b-\varphi\circ \gamma)^2} + \frac{2|B|^2\left(\D -\p_t\right)(\varphi\circ \gamma)}{(b-\varphi\circ \gamma)^3} \\
& + \frac{2\la \n f, \n(\varphi\circ \gamma) \ra}{b-\varphi\circ \gamma} + \frac{2|B|^2|\n(\varphi\circ \gamma)|^2}{(b-\varphi\circ \gamma)^4} \\
& +  \frac{2\la \n |B|^2, \n(\varphi\circ \gamma) \ra}{(b-\varphi\circ \gamma)^3}.
\endaligned
\end{equation}
The inequality  (\ref{eqn-Para1phi})and Theorem A in \cite{Wan03} imply 
\begin{equation}\label{eqn-Pa4}\aligned
\left(\D - \p_t \right)(\varphi\circ \gamma) = &\sum_{i=1}^n {\rm Hess} (\varphi) (d\gamma(e_i), d\gamma(e_i))\circ \gamma 
+ d\varphi(\tau(\gamma)-\p_t \gamma) \\
=& \sum_{i=1}^n {\rm Hess} (\varphi) (d\gamma(e_i), d\gamma(e_i)) \circ \gamma \\
= & \cos(\rho\circ \gamma)|B|^2.  
\endaligned
\end{equation}
By Corollary 3.5 in \cite{Hui84}, we have
\begin{equation}\label{eqn-Pa5}\aligned
\left(\D - \p_t \right)|B|^2 = 2|\n B|^2 - 2|B|^4 \geq &  2|\n |B||^2 - 2|B|^4.
\endaligned
\end{equation}
Substituting (\ref{eqn-Pa4}) and (\ref{eqn-Pa5}) into (\ref{eqn-Pa3}), we have
\begin{equation}\label{eqn-Pa6}\aligned
\left(\D - \p_t \right) f \geq & 2(b - \varphi\circ \gamma)^2 \left( \frac{\cos(\rho\circ \gamma)}{b-\varphi\circ \gamma}-1 \right) f^2 + \frac{2|\n |B||^2}{(b-\varphi\circ \gamma)^2} \\
&+ \frac{2\la \n f, \n(\varphi\circ \gamma) \ra}{b-\varphi\circ \gamma} 
+ \frac{2|B|^2|\n(\varphi\circ \gamma)|^2}{(b-\varphi\circ \gamma)^4} +  \frac{2\la \n |B|^2, \n(\varphi\circ \gamma) \ra}{(b-\varphi\circ \gamma)^3}.
\endaligned
\end{equation}
By the Schwartz inequality, we obtain
\begin{equation}\label{eqn-Pa7}\aligned
\left| \frac{2\la \n |B|^2, \n(\varphi\circ \gamma) \ra}{(b-\varphi\circ  \gamma)^3} \right| \leq \frac{2|\n |B||^2}{(b-\varphi\circ \gamma)^2} + \frac{2|B|^2|\n(\varphi\circ \gamma)|^2}{(b-\varphi\circ \gamma)^4}.
\endaligned
\end{equation}
Combining (\ref{eqn-Pa6}) with (\ref{eqn-Pa7}),  it follows
\begin{equation}\label{eqn-Pa8}\aligned
\left(\D - \p_t \right) f \geq & 2(1-b)(b - \varphi\circ \gamma)f^2 
+ \frac{2\la \n f, \n(\varphi\circ \gamma) \ra}{b-\varphi\circ \gamma}. 
\endaligned
\end{equation}

Let $\eta(r, t): \mathbb{R} \times \mathbb{R} \to \mathbb{R}$ be a smooth function supported on $[-R, R]\times [-T, 0]$, satisfying the following properties:
\begin{itemize}
    \item[(1)] $\eta(r, t) \equiv 1$ on $[-\frac{R}{2}, \frac{R}{2}] \times [-\frac{T}{2}, 0]$ and $0\leq \eta \leq 1.$
    \item[(2)] $\eta(r, t)$ is decreasing if $r\geq 0$, i.e., $\p_r \eta \leq 0$.
    \item[(3)] $\frac{|\p_r \eta|}{\eta^a}\leq \frac{C_a}{R}, \frac{|\p_r^{2}\eta|}{\eta^a} \leq \frac{C_a}{R^2}$ for $a=\frac{1}{2}, \frac{3}{4}$.
    \item[(4)] $\frac{|\p_t \eta|}{\eta^{\frac{1}{2}}} \leq \frac{C}{T}$.
\end{itemize}
Such a function was explicitly constructed in \cite{Kun16} (see also \cite{LY, SZ}). 
Let $\phi:= \eta(r(F), t)$, where $r(F):= |F|$.

Let $L:= -\frac{2\n  (\varphi\circ \gamma)}{ b - \varphi\circ \gamma}$. Then we get
\begin{equation}\label{eqn-Pa9}\aligned
&\left(\D- \p_t \right)(\phi f) + \la L, \n(\phi f) \ra -2 \left\la \frac{\n \phi}{\phi}, \n(\phi f) \right\ra \\
=& f(\D-\p_t)\phi + \phi (\D-\p_t) f + 2\la \n \phi, \n f \ra \\
&+ \phi\la L, \n f \ra + f\la L, \n \phi \ra - \frac{2|\n \phi|^2 f}{\phi} - 2\la \n \phi, \n f \ra \\
=& f(\D-\p_t)\phi + \phi (\D-\p_t) f -\frac{2\phi \la \n f, \n (\varphi\circ \gamma)\ra}{ b -\varphi \circ \gamma } \\
& - \frac{2f \la \n \phi, \n (\varphi \circ \gamma) \ra}{ b-\varphi\circ \gamma } - \frac{2|\n \phi|^2f}{\phi}.
\endaligned
\end{equation}
From (\ref{eqn-Pa8}) and (\ref{eqn-Pa9}), we obtain
\begin{equation}\label{eqn-Pa10}\aligned
&\left(\D- \p_t \right)(\phi f) + \la L, \n(\phi f) \ra -2 \left\la \frac{\n \phi}{\phi}, \n(\phi f) \right\ra \\
\geq & 2(1-b)(b - \varphi\circ \gamma)\phi f^2 + f \left( \D-\p_t \right)\phi \\
& - \frac{2f \la \n \phi, \n(\varphi\circ \gamma) \ra}{b-\varphi\circ \gamma} - \frac{2|\n \phi|^2f}{\phi}.
\endaligned
\end{equation}

By Proposition 3.3 in \cite{Kun15}, $M_t$ can be written as a complete graph, thus $D_{R,T}(o)$ is compact.  
Hence $\phi f$ attains its maximum at some point $(p_1, t_1)$ in $D_{R, T}(o)$. By the maximum principle, at $(p_1, t_1)$, we have
\[
\n(\phi f) = 0, \quad \D(\phi f) \leq 0, \quad \p_t(\phi f) \geq 0.
\]
Then from (\ref{eqn-Pa10}), we derive at $(p_1, t_1)$ ,
\begin{equation}\label{eqn-para8}\aligned
2(1-b)(b - \varphi\circ \gamma)\phi f^2 \leq  - f(\D-\p_t)\phi 
+ \frac{2f \la \n \phi, \n(\varphi\circ \gamma) \ra}{b-\varphi\circ \gamma} + \frac{2|\n \phi|^2f}{\phi}.
\endaligned
\end{equation}

Direct computation gives
\begin{equation*}\label{eqn-main}\aligned
(\D-\p_t) r^2 = (\D-\p_t) |F|^2 = 2n.
\endaligned
\end{equation*}
and
\[
(\D-\p_t) r^2 = 2|\n r|^2 +2r(\D-\p_t) r \geq 2r(\D-\p_t )r,
\]
The above two formulas imply that 
\[
(\D-\p_t)r \leq \frac{n}{r}.
\]
It follows that
\begin{equation}\label{eqn-Para9}\aligned
 - f(\D-\p_t)\phi = & -f\p_r^{2}\eta |\n r|^2 - f\p_r\eta (\D-\p_t) r + f \p_t\eta \\
 \leq & -f\p_r^{2}\eta |\n r|^2 - f\p_r\eta \cdot\frac{n}{r} + f \p_t\eta \\
 \leq & f|\p_{r}^2 \eta| + f |\p_r \eta|\cdot \frac{2n}{R} + f|\p_t \eta| \\
 =& f\phi^{\frac{1}{2}} \cdot \frac{|\p_r^{2}\eta|}{\phi^{\frac{1}{2}}} + f\phi^{\frac{1}{2}} \cdot\frac{|\p_r \eta| \frac{2n}{R}}{\phi^{\frac{1}{2}}} +  f\phi^{\frac{1}{2}} \cdot \frac{|\p_t\eta|}{\phi^{\frac{1}{2}}} \\
 \leq &  \frac{\e}{5} \phi f^2 + \frac{5}{4\e} \left( \frac{|\p_r^{2}\eta|}{\phi^{\frac{1}{2}}} \right)^2 +  \frac{\e}{5} \phi f^2 +\frac{5}{4\e} \left(\frac{2n}{R}\right)^2 \cdot\left(\frac{ |\p _r \eta|}{\phi^{\frac{1}{2}}} \right)^2 \\
 &+  \frac{\e}{5} \phi f^2 +\frac{5}{4\e} \left( \frac{|\p_t\eta|}{\phi^{\frac{1}{2}}} \right)^2 \\
 \leq &  \frac{3\e}{5} \phi f^2 + \frac{20n^2+5}{4\e} \frac{C}{R^4} +  \frac{5}{4\e} \frac{C}{T^2}.
\endaligned
\end{equation}
(Note that since $\p_r\eta=0$ for $r\leq \frac{R}{2}$, then we may assume that $r\geq \frac{R}{2}$ in the second inequality).

The Schwartz inequality implies 
\begin{equation}\label{eqn-Pa11}\aligned
\frac{2f \la \n \phi, \n(\varphi\circ \gamma) \ra}{b-\varphi\circ \gamma} \leq & \frac{2f|\n \phi|\cdot |B|}{b-\varphi\circ \gamma} =2f^{\frac{3}{2}}|\n \phi| \\
 = &2f^{\frac{3}{2}}\phi^{\frac{3}{4}} \cdot \frac{|\n \phi|}{\phi^{\frac{3}{4}}} \leq \frac{\e}{5}\phi f^2 +\frac{3375}{16\e^3} \frac{C}{R^4}.
\endaligned
\end{equation}

The last term of the right hand side of (\ref{eqn-para8}) can be estimated as follows
\begin{equation}\label{eqn-Para11}\aligned
\frac{2|\n \phi|^2f}{\phi} = f\phi^{\frac{1}{2}}\cdot \frac{2|\n \phi|^2}{\phi^{\frac{3}{2}}} \leq \frac{\e}{5} \phi f^2 +\frac{5}{\e} \left( \frac{|\n\phi|}{\phi^{\frac{3}{4}}} \right)^4 \leq \frac{\e}{5} \phi f^2 +\frac{5}{\e} \frac{C}{R^4}.
\endaligned
\end{equation}
Substituting (\ref{eqn-Para9}),  (\ref{eqn-Pa11}) and (\ref{eqn-Para11}) into (\ref{eqn-para8}), we have
\begin{equation*}\label{eqn-Para12}\aligned
2(1-b)(b - \varphi\circ \gamma(p_1,  t_1)) \phi f^2 \leq  \e \phi f^2 +\left( \frac{20n^2+25}{4\e}+ \frac{3375}{16\e^3} \right) \frac{C}{R^4} + \frac{5}{4\e}\frac{C}{T^2}.
\endaligned
\end{equation*}

Choosing 
\[
\e=(1-b)(b - \varphi\circ \gamma(p_1,  t_1)).
\]
Then we get at $(p_1, t_1)$,
\begin{equation}\label{eqn-Pa12}\aligned
 \phi f^2 \leq  C \left( \left( \frac{1}{\e^2}+ \frac{1}{\e^4}  \right) \frac{1}{R^4} + \frac{1}{\e^2} \frac{1}{T^2} \right),
\endaligned
\end{equation}
where $C$ is a positive constant which is independent of $R$ and $T$.

Since 
\begin{equation*}\label{eqn-QH0}
\lim_{t\to (\frac{\pi}{2})^-}\frac{(\cos t)^{-1}}{\left( \frac{\pi}{2} -t \right)^{-1}} = \lim_{t\to (\frac{\pi}{2})^-}\frac{\frac{\pi}{2} -t}{\cos t} =  \lim_{t\to (\frac{\pi}{2})^-} \frac{-1}{-\sin t} = 1.
\end{equation*}
Therefore we obtain
\[
\frac{1}{\e}=\frac{1}{(1-b)(b - \varphi\circ \gamma(p_1,  t_1)) } \leq 4 \sup_{D_{R, T}(o)} \left(\cos(\rho\circ \gamma)\right)^{-2} \leq C\sup_{D_{R, T}(o)} \left( \frac{\pi}{2} - \rho\circ \gamma \right)^{-2}.
\]
Hence from (\ref{eqn-Pa12}),  we derive
\begin{equation*}\label{eqn-Pa131}
\aligned
 \sup_{D_{\frac{R}{2}, \frac{T}{2}}(o)} \frac{|B|}{b - \varphi\circ \gamma} 
\leq & C \left( \frac{1}{R}\cdot \left( \sup_{D_{R, T}(o)}\left( \frac{\pi}{2} - \rho\circ \gamma \right)^{-1} \right.\right. \\
& \left.\left.+  \sup_{D_{R, T}(o)}\left( \frac{\pi}{2} - \rho\circ \gamma \right)^{-2} \right) 
+\frac{1}{\sqrt T}\cdot  \sup_{D_{R, T}(o)}\left( \frac{\pi}{2} - \rho\circ \gamma \right)^{-1} \right).
\endaligned
\end{equation*} 
\qed

\noindent{\bf Proof of Theorem \ref{thm-main2}} By (\ref{eqn-Pa01}), we have
\begin{equation}\label{eqn-Pa13}
\aligned
\frac{1}{b}\sup_{D_{\frac{R}{2}, \frac{R^2}{2}}(o)} |B|
\leq  C \left( \frac{1}{R}\left(\sup_{D_{R, R^2}(o)}  \left( \frac{\pi}{2} - \rho\circ u \right)^{-1}
+ \sup_{D_{R, R^2}(o)}  \left( \frac{\pi}{2} - \rho\circ u \right)^{-2} \right)\right).
\endaligned
\end{equation}
Letting $R \to +\infty$ in (\ref{eqn-Pa13}), then we derive $B \equiv 0$.  Namely,  $M_t$ has to be affine linear for any $t\in (-\infty, 0]$. 

\vskip24pt

 \noindent{\bf Final Remark} Theorem \ref{thm-main2} also holds true for eternal solutions to the mean curvature flows if we take the time interval as $[-T, T]$ in the proof of Theorem \ref{thm-main1}.
 
 \vskip12pt

  \noindent{\bf Acknowledgements}  This work is partially supported by NSFC (Nos. 12331002, 12471050) and Hubei Provincial Natural Science Foundation of China (No. 2024AFB746).  The authors would like to thank Professors Guofang Wang and Yong Luo for helpful discussions.  

\vskip24pt

\end{document}